\documentclass[12pt]{amsart}
\usepackage{amsfonts,amsmath,amssymb,amsthm}
\usepackage{latexsym}
\usepackage{mathrsfs}
\usepackage{tikz}
\usepackage{color}
\usepackage{cite}

\newtheorem{prop}{Proposition}[section]
\newtheorem{thm}{Theorem}[section]
\newtheorem{rem}{Remark}
\newtheorem{corollary}{Corollary}[section]

\newtheorem{conj}{Conjecture}[section]

\newcommand{\half}{{\textstyle{\frac{1}{2}}}}

\begin{document}

\title[Brauer and Bannai--Ito algebras]{Centralizers of the superalgebra $\mathfrak{osp}(1|2)$: the Brauer \vspace{2mm}\\algebra as a quotient of the Bannai--Ito algebra}

\author[N.Cramp\'e]{Nicolas Cramp\'e$^{\dagger,*}$}
\address{$^\dagger$ Institut Denis-Poisson CNRS/UMR 7013 - Universit\'e de Tours - Universit\'e d'Orl\'eans, 
Parc de Grammont, 37200 Tours, France.}
\email{crampe1977@gmail.com}

\author[L.Frappat]{Luc Frappat$^{\ddagger}$}
\address{$^\ddagger$ Laboratoire d'Annecy-le-Vieux de Physique Th\'eorique LAPTh, Univ. Grenoble Alpes, 
Univ. Savoie Mont Blanc, CNRS, F-74000 Annecy, France. }
\email{luc.frappat@lapth.cnrs.fr}

\author[L.Vinet]{Luc Vinet$^{*}$}
\address{$^*$ Centre de recherches math\'ematiques, Universit\'e de Montr\'eal,
P.O. Box 6128, Centre-ville Station,
Montr\'eal (Qu\'ebec), H3C 3J7, Canada.}
\email{vinet@CRM.UMontreal.ca}

\begin{abstract}
We provide an explicit isomorphism between a quotient of the Bannai--Ito algebra and the Brauer algebra.
We clarify also the connection with the action of the Lie superalgebra $\mathfrak{osp}(1|2)$ on the threefold
tensor product of its fundamental representation. Finally, a conjecture is proposed to describe the centralizer 
of $\mathfrak{osp}(1|2)$ acting on three copies of an arbitrary finite irreducible representation in terms of a quotient of the Bannai--Ito algebra.
\end{abstract}

\maketitle

\textit{To the fond memory of Peter Freund, a much esteemed scientist who always generously shared his immense culture.}

\vspace{3mm}

\section{Introduction}

The purpose of this paper is to obtain the relation between the Bannai--Ito  and  the Brauer algebras.
The Brauer algebra has been introduced in \cite{Brauer} in the framework of the Schur--Weyl duality for the orthogonal and symplectic groups whereas the Bannai--Ito algebra has been defined in \cite{TVZ}
to give an algebraic description of the eponym polynomials \cite{BI}.

Both algebras are associated to the centralizer of the Lie superalgebra $\mathfrak{osp}(1|2)$; a connection between these algebras is therefore expected. 
Inspired by recent results relating the Racah algebra and the centralizers of $\mathfrak{su}(2)$ \cite{CPV}, we have found an 
explicit isomorphism (that we report here) between a quotient of the Bannai--Ito algebra, the Brauer algebra and the action of $\mathfrak{osp}(1|2)$ on the threefold tensor product
of fundamental representations of this superalgebra. The quotient of the Bannai--Ito algebra is linked to the direct sum decomposition of the tensor product 
of the three fundamental representations which is usefully depicted with the help of the associated Bratteli diagram. This construction can be generalized so as to give a conjecture 
which if true, describes in terms of generators and relations the centralizers of the action of $\mathfrak{osp}(1|2)$ on the tensor product
of three arbitrary irreducible representations of finite dimension.

The outline of the paper is as follows. In Section \ref{sec:BIB}, we define the universal Bannai--Ito and Brauer algebras. 
Then, we state and prove the main theorem of this paper which gives the isomorphism between a quotient of the Bannai--Ito algebra and a specialization of the Brauer algebra.
Section \ref{sec:BIosp} describes the relation that the previous result has with the centralizer of the threefold tensor 
product of the fundamental representation of the Lie superalgebra $\mathfrak{osp}(1|2)$. 
We first recall the definition and the algebraic properties of $\mathfrak{osp}(1|2)$ in Subsection \ref{sec:osp}. 
We then briefly give an overview of the finite irreducible representations  of $\mathfrak{osp}(1|2)$ in Subsection \ref{sec:irreps}.
The connection between the centralizer of three fundamental representations of $\mathfrak{osp}(1|2)$ and the Bannai--Ito algebra is explained in 
Subsection \ref{sec:centr}. Finally, we propose in Section \ref{sec:conj}  a conjecture for an isomorphism between 
the centralizer for the threefold tensor product of an arbitrary  $\mathfrak{osp}(1|2)$ finite irreducible representation with itself and a quotient of the Bannai--Ito algebra. 

\section{Quotient of the Bannai--Ito algebra and Brauer algebra \label{sec:BIB}}

The universal Bannai--Ito algebra  $I_3$  is generated by three generators $X$, $Y$ and $Z$ and three central elements $\omega_X$, $\omega_Y$ and $\omega_Z$ satisfying 
the relations \cite{BI,TVZ,GVZa,DGTVZ,DGV2}
\begin{subequations}\label{eq:BI}
\begin{align}
& \{X,Y\}=X+Y+Z+\omega_Z \,, \label{eq:BI1} \\ 
& \{X,Z\}=X+Y+Z+\omega_Y \,, \label{eq:BI2} \\
& \{Y,Z\}=X+Y+Z+\omega_X \,, \label{eq:BI3}
\end{align}
\end{subequations}
where $\{X,Y\}=XY+YX$. The usual presentation of the Bannai--Ito algebra \cite{GVZa,DGTVZ,DGV2} is easily recovered by the following affine transformations: 
$X\rightarrow X+1/2$, $Y\rightarrow Y+1/2$, $Z\rightarrow Z+1/2$, $\omega_X\rightarrow \omega_X-1$, $\omega_Y\rightarrow \omega_Y-1$ and $\omega_Z\rightarrow \omega_Z-1$.

The Brauer algebra $B_3(\eta)$ is the unital algebra generated by the four generators $s_1$, $s_2$, $e_1$ and $e_2$ with the defining relations \cite{Brauer}
\begin{alignat}{3}
& s_i^2=1\,, &\qquad  & e_i^2=\eta e_i\,, &\qquad  & s_ie_i=e_is_i=e_i\,,\\
& s_1s_2s_1=s_2s_1s_2\,, &\qquad  & e_1e_2e_1=e_1\,, &\qquad  & e_2e_1e_2=e_2\,,\\
& s_1e_2e_1=s_2e_1\,, &\qquad  & e_2e_1s_2=e_2s_1\,. &\qquad  &
\end{alignat}
Let us recall that the dimension of the Brauer algebra $B_3(\eta)$ is $15$ and it is easy to prove that the following relations also hold:
\begin{alignat}{4}
& s_1s_2e_1=e_2e_1\,, &\qquad  & e_2s_1s_2=e_2e_1\,, &\qquad  & s_2e_1s_2=s_1e_2s_1\,, &\qquad  & \\
& s_2e_1e_2=s_1e_2\,, &\qquad  & e_1e_2s_1=e_1s_2\,, &\qquad  & e_1s_2e_1=e_1\,, &\;\;  & e_2s_1e_2=e_2\,, \\
& s_2s_1e_2=e_1e_2\,, &\qquad  & e_1s_2s_1=e_1e_2\,. &\qquad  & &\qquad  & 
\end{alignat}

The main result of this paper is stated in the following theorem where we give an explicit isomorphism between a quotient of the Bannai--Ito algebra and a specialization of the Brauer algebra.
\begin{thm}\label{thm:iso}
The quotient of the Bannai--Ito algebra $I_3$ by the following relations $\omega_X=\omega_Y=\omega_Z=\omega$ and
\begin{align}
& X\left(X^2-4\right)=0\ , \qquad Y\left(Y^2-4\right)=0 \ , \qquad  Z\left(Z^2-4\right)=0 \ , \label{eq:q1}\\
&\left(\omega+4\right)\left(\omega-2\right)\left(\omega-8\right)\left(\omega-14\right)=0\ , \label{eq:q2} \\ 
& \left(X-\omega+16\right) \left(X-\omega+10\right) \left(X-\omega+8\right) \left(X-\omega+6\right) \left(X-\omega\right) \left(X-\omega-2\right) \left(X-\omega-6\right)=0,\label{eq:q3}\\
& \left(Y-\omega+16\right) \left(Y-\omega+10\right) \left(Y-\omega+8\right) \left(Y-\omega+6\right) \left(Y-\omega\right) \left(Y-\omega-2\right) \left(Y-\omega-6\right)=0,\label{eq:q4}\\
& \left(Z-\omega+16\right) \left(Z-\omega+10\right) \left(Z-\omega+8\right) \left(Z-\omega+6\right) \left(Z-\omega\right) \left(Z-\omega-2\right) \left(Z-\omega-6\right)=0.\label{eq:q5}
\end{align}
denoted  $\overline{I}_3$,
is isomorphic to the Brauer algebra $B_3(-1)$.
This isomorphism is given explicitly by:
\begin{align}
 \Psi\  :\ \overline{I}_3&\rightarrow B_3(-1)\\
 X&\mapsto 2(s_1+e_1)\nonumber\\
 Y&\mapsto 2(s_2+e_2)\nonumber\\
 Z&\mapsto 2s_2(s_1+e_1)s_2=2 s_1(s_2+e_2)s_1.\nonumber
\end{align}
The image of $\omega$ by $\Psi$ is given by $\Psi(\omega)=\{\Psi(X),\Psi(Y)\}-\Psi(X)-\Psi(Y)-\Psi(Z)$.
\end{thm}

\proof The first step of the proof consists in proving that $\Psi$ is a homomorphism, in other words that $\Psi(X)$, $\Psi(Y)$, $\Psi(Z)$ and $\Psi(\omega)$ satisfy the relations of the quotient $\overline{I}_3$ of the Bannai--Ito algebra.
The relations \eqref{eq:BI1} and \eqref{eq:q1} are easy to check.
The relation \eqref{eq:BI2} gives
\begin{align}
 \{\Psi(X)&,\Psi(Z)\}-\Psi(X)-\Psi(Y)-\Psi(Z)-\Psi(\omega)=4\{s_1+e_1,s_1(s_2+e_2)s_1\}-4\{s_1+e_1,s_2+e_2\}\nonumber\\
 &=4\left(  (s_2+e_2)s_1+s_1(s_2+e_2)+\underbrace{e_1(s_2+e_2)s_1}_{=e_1e_2+e_1s_2}+\underbrace{s_1(s_2+e_2)e_1}_{=e_2e_1+s_2e_1} -\{s_1+e_1,s_2+e_2\}  \right)=0 \,.
\end{align}
The relation \eqref{eq:BI3} is computed similarly. Relations \eqref{eq:q2} and \eqref{eq:q5} need more work to be verified. 
We prove them in the faithful $15\times 15$ regular representation of the Brauer algebra.

The second step is to show that $\Psi$ is surjective which is done easily. Indeed, one gets $\Psi(1+X/2-X^2/4)=s_1$, $\Psi(1+Y/2-Y^2/4)=s_2$, $\Psi(-1+X^2/4)=e_1$ and $\Psi(-1+Y^2/4)=e_2$.
The 4 generators of the Brauer algebra belong to the image of $\Psi$. Therefore, $\Psi$ is surjective.

The last step requires demonstrating that $\Psi$ is injective. We know that the dimension of $B_3(-1)$ is 15. To prove the injectivity, we have to show that there is a generating family of generators of dimension 15 in $\overline I_3$. 
By using the anti-commutation relations \eqref{eq:BI}, \eqref{eq:q1} and \eqref{eq:q2}, it is easy to see that the following ensemble
\begin{equation}
 \{ w^j X^x Y^y Z^z \ | \ j=0,1,2,3\text{ and } x,y,z=0,1,2 \} \label{eq:genf1}
\end{equation}
forms a set of generators. We will show that there exist supplementary relations between the elements of that set.
Compute $X^2 \eqref{eq:BI1} - X \eqref{eq:BI1} X + \eqref{eq:BI1} X^2$, using the fact that $X^3=4X$ in $\overline I_3$, it is seen that the following relation is implied 
in $\overline I_3$
\begin{equation}
X^2Z=-X^2Y-\frac{1}{3}(\omega-2) (X^2-2X)+2XY+2XZ\;.\label{eq:res1}
\end{equation}
Similarly, one gets
\begin{eqnarray}
 YZ^2&=&-XZ^2-\frac{1}{3}(\omega-2) (Z^2-2Z)+2 XZ+2 Y Z\ ,\label{eq:res2}\\
 Y^2Z&=&-XY^2-\frac{1}{3}(\omega-2) (Y^2-2Y)+2\omega +2X+2Y+2Z\ .\label{eq:res3}
\end{eqnarray}
The equations \eqref{eq:res4} to \eqref{eq:res8} below are obtained as follows:
multiplying expression \eqref{eq:res1} by $Y$ on the right and ordering with \eqref{eq:BI}, one finds \eqref{eq:res4}; using this last equation after having multiplied \eqref{eq:res1} by $Z$ on the right leads to \eqref{eq:res5}; multiplying  \eqref{eq:res2} by $X$ on the left and using  \eqref{eq:res5}, one arrives at \eqref{eq:res6}; multiplying  \eqref{eq:res2} by $Y$ on the left and calling upon \eqref{eq:res6}, one gets \eqref{eq:res7} and finally, one obtains \eqref{eq:res8} by multiplying \eqref{eq:res3} by $X$ on the left.
\begin{eqnarray}
 X^2YZ&= &X^2Y^2  +\frac{1}{3} (\omega-2)( X^2Y-2XY+2X^2-4X) -2XY^2+2XYZ\ , \label{eq:res4} \\
 X^2Z^2&=& -X^2Y^2 +\frac{1}{9}(\omega-8)(\omega-2)(X^2-2X)+2XY^2+2XZ^2\ ,\label{eq:res5} \\
 XYZ^2&=&X^2Y^2-\frac{1}{9}(\omega-2)^2(X^2 -2X) -\frac{1}{3}(\omega+4)(XZ^2-2XZ)\nonumber\\
 &&-2XY^2+4XY+2XYZ-2X^2Y \ , \label{eq:res6} \\
 Y^2Z^2&=& X^2Y^2-\frac{1}{9}(\omega-2)^2(X^2 -2X)+\frac{1}{9}(\omega-8)(\omega-2)(Z^2-2Z)-4XY^2+4XY\nonumber\\
 &&-\frac{2}{3}(\omega-2)(Y^2-2Y)-2X^2Y-2XZ^2+4XZ+4X+4Y+4Z+4\omega\ ,\label{eq:res7} \\
 XY^2Z&=&-X^2Y^2-\frac{1}{3}(\omega-2) (XY^2-2XY)+2X^2+2XY+2XZ+2\omega X\ . \label{eq:res8} 
\end{eqnarray}
Multiplying \eqref{eq:res8} by $X$ on the left and using \eqref{eq:res1}, 
multiplying \eqref{eq:res6} by $X$ on the left and using \eqref{eq:res1}, \eqref{eq:res4}, \eqref{eq:res5}, 
multiplying \eqref{eq:res7} by $X$ on the left and using \eqref{eq:res1}, \eqref{eq:res5}, 
one finds in a similar way
\begin{eqnarray}
 X^2Y^2Z&=&-4XY^2-\frac{1}{3}(\omega-2) (X^2Y^2-2X^2Y+2X^2-4X)\nonumber\\
     &&+4XY+4XZ+2\omega X^2+8X\ , \label{eq:res9} \\
 X^2YZ^2&=&4XYZ -\frac{1}{27}(\omega+4)(\omega-2)(\omega-8)(X^2-2X)\nonumber\\
 &&+\frac{1}{3}(\omega+4)(X^2Y^2-2XY^2+2XZ^2+4XZ)\ ,\label{eq:res10}\\
 XY^2Z^2&=&-2X^2Y^2+\frac{1}{9}(\omega+4)(\omega-14)(XZ^2-2XZ)-\frac{2}{3}(\omega+4)(XY^2-2XY)\nonumber\\
 &&+4(XZ^2-XZ+XY^2-XY+X^2+\omega X)\ .\label{eq:res11}
\end{eqnarray}
One also obtains
\begin{eqnarray}
&& X^2Y^2Z^2=\frac{2(\omega+4)(\omega-14)}{9}(XY^2-2XY+XZ^2-2XZ   ) +8XZ^2-8XY-8XZ  \nonumber\\
 &&+\frac{(\omega-8)^2(\omega-2)(\omega+4)}{81}(X^2-2X)-\frac{(\omega+4)(\omega-8)}{9}(X^2Y^2-2X^2Y)+8X^2 +8\omega X  \ . \label{eq:res12}
\end{eqnarray}
Then, we deduce from relations \eqref{eq:res1} to \eqref{eq:res12} that the generating family \eqref{eq:genf1} of $\overline I_3$ reduces to 
\begin{equation}
 \mathcal{C}=\mathcal{F} \cup \omega \mathcal{F}\cup \omega^2 \mathcal{F}\cup \omega^3 \mathcal{F}\ ,
\end{equation}
where
\begin{equation}
 \mathcal{F}=\{1 ,\ X ,\ Y,\ Z,\ X^2,\ Y^2,\ Z^2,\ XY,\ XZ ,\ YZ ,\ X^2Y ,\ XY^2,\ XZ^2 ,\ XYZ,\ X^2Y^2\}\ .
\end{equation}
To find supplementary relations between the 60 elements of the set $\mathcal{C}$, we construct the regular action of the generators $X$, $Y$, $Z$ and $\omega$ on $\mathcal{C}$
thereby associating to each of these 4 generators a $60\times 60$ matrix. Demanding that these matrices satisfy the relations of the quotiented Bannai--Ito algebra, 
we find $45$ constraints.
An abstract mathematical software has been useful to perform these computations. 
We thus deduce that $\mathcal{F}$ is a generating family. Since $\mathcal{F}$ has 15 elements, this implies the injectivity of 
$\Psi$ and concludes the proof.
\endproof

\begin{rem}
 From the previous theorem, we know that the dimension of $\overline I_3$ is $15$. With computations similar to those used in the proof, we can also show that 
 the dimensions of the quotients of $\overline I_3$ by $\omega=-4$, 
  $\omega=2$, 
  $\omega=8$ 
 and 
  $\omega=14$ 
 are respectively $4$, $1$, $9$ and $1$. 
\end{rem}

While the relations \eqref{eq:q1}--\eqref{eq:q5} used to define the quotient of the Bannai--Ito algebra seem artificial at first glance, we are going to show in the following that this quotient is natural when we consider the image of the Bannai--Ito algebra in three copies 
of the fundamental representation of the Lie superalgebra $\mathfrak{osp}(1|2)$.

\section{Bannai--Ito algebra and Lie superalgebra $\mathfrak{osp}(1|2)$ \label{sec:BIosp}}

\subsection{Algebraic definitions and properties \label{sec:osp}}

In this subsection, we recall definitions and well-known results concerning the Lie superalgebra $\mathfrak{osp}(1|2)$.

This superalgebra has two odd generators $F^\pm$ and three even generators $H$, $E^\pm$ satisfying the following 
(anti-)commutation relations \cite{kac}
\begin{alignat}{2}
& [H,E^\pm]=\pm E^\pm \,,             &\quad & [E^+,E^-]=2H\,,\\
& [H,F^\pm]=\pm \frac{1}{2} F^\pm \,, &\quad & \{F^+,F^-\}=\frac{1}{2}H\,,\\
& [E^\pm,F^\mp]=-F^\pm \,,            &\quad & \{F^\pm,F^\pm\}=\pm\frac{1}{2} E^\pm\;.
\end{alignat}
Remark that the subalgebra generated by $H$ and $E^\pm$ is isomorphic to $\mathfrak{su}(2)$.
The $\mathbb{Z}_2$-grading of $\mathfrak{osp}(1|2)$ can be encoded by the grading involution $R$ satisfying
\begin{equation}
 [R,E^\pm]=0 \;, \qquad  [R,H]=0 \;, \qquad  \{R,F^\pm\}=0 \qquad \text{and} \qquad  R^2=1\;.
\end{equation}
In the universal enveloping algebra $U(\mathfrak{osp}(1|2))$, one defines the sCasimir by \cite{Pinc,Les}
\begin{equation}
 S=[F^+,F^-]+\frac{1}{8}\;.
\end{equation}
It anti-commutes with the odd generators, $\{S,F^\pm\}=0$ and commutes with the even ones, $[S,E^\pm ]=0$, $[S,H]=0$.
A central element $Q$ of $U(\mathfrak{osp}(1|2))$ can be constructed as follows by using the sCasimir and the grading involution:
\begin{equation}
 Q=S\;R=[F^+,F^-]R+\frac{R}{8}\;.\label{eq:Q}
\end{equation}
The $U(\mathfrak{osp}(1|2))$ algebra is endowed with a coproduct $\Delta$ defined as the algebra homomorphism 
satisfying
\begin{alignat}{2}
& \Delta(E^\pm)=E^\pm \otimes 1 +1 \otimes E^\pm \,, &\qquad &\Delta(H)=H\otimes 1 + 1\otimes H\,, \\ 
& \Delta(F^\pm)=F^\pm \otimes R +1 \otimes F^\pm \,, &\qquad &\Delta(R)=R\otimes R\;.
\end{alignat}
Denote by $U^3$ the threefold tensor product $U(\mathfrak{osp}(1|2))\otimes U(\mathfrak{osp}(1|2)) \otimes  U(\mathfrak{osp}(1|2))$. We define in $U^3$ the following algebraic elements
\begin{alignat}{3}
 &Q_1=Q\otimes 1 \otimes 1 \,, &\qquad &Q_2=1\otimes Q \otimes 1 \,, &\qquad &Q_3=1\otimes 1 \otimes Q \,,\\
 &Q_{12}=\Delta(Q) \otimes 1 \,, &\qquad &Q_{23}=1 \otimes \Delta(Q)\,, &\qquad & \\
 &Q_4=(\Delta \otimes 1)\Delta(Q)\,.
\end{alignat}
Finally, one introduces also
\begin{equation}
 Q_{13}=\left([F^+\otimes R \otimes R+1\otimes 1 \otimes F^+,F^-\otimes R \otimes R+1\otimes 1 \otimes F^-]+\frac{1}{8}\right)R\otimes 1 \otimes R\;.
\end{equation}
The relations between the Bannai--Ito algebra and the algebraic elements $Q$ are given by the following statement \cite{GVZa,DGTVZ}:
\begin{prop}\label{pr:Phi}
The map $\Phi : I_3 \to U^3$ defined by 
\begin{equation}
X \mapsto -4Q_{12}+\frac{1}{2}\,,\qquad  Y \mapsto -4Q_{23}+\frac{1}{2}\,,\qquad Z \mapsto -4Q_{13}+\frac{1}{2}\,,\label{eq:Phi1}
\end{equation}
and
\begin{subequations}
\begin{align}
& \omega_X \mapsto 32(Q_1Q_2+Q_3Q_4)-1 \,, \\
& \omega_Y \mapsto 32(Q_2Q_3+Q_1Q_4)-1 \,, \\
& \omega_Z \mapsto 32(Q_1Q_3+Q_2Q_4)-1\;,
\end{align}
\label{eq:Phi2}
\end{subequations}
is an algebra homomorphism. 
\end{prop}
Note that the shift by 1/2 in \eqref{eq:Phi1} is due to our definition of the $X,Y,Z$ generators in comparison to \cite{GVZa,DGTVZ}. 
The importance of the previous construction lies in the fact that the image of the Bannai--Ito algebra by $\Phi$ 
belongs to the centralizer of the image of $U(\mathfrak{osp}(1|2))$ by $(\Delta \otimes 1)\Delta$.
Indeed, one gets
\begin{equation}\label{eq:centu2}
 [ (\Delta \otimes 1)\Delta(g) ,  \Phi(x) ] =0 \qquad  \forall g\in  U(\mathfrak{osp}(1|2)) \quad\text{and}\quad \forall x\in I_3\;.
\end{equation}

\subsection{Finite irreducible representations of $\mathfrak{osp}(1|2)$ \label{sec:irreps}} 

The finite irreducible representations $[j]^\pm$ of $\mathfrak{osp}(1|2)$ are labeled 
by an integer or an half integer $j$ but also by a sign $\pm$ corresponding to the parity of the highest weight 
($+$ stands for a bosonic state and $-$ for the fermionic state) \cite{PA1975,SNR1977,ERSS}.
More precisely, if we denote by $v_j^\pm$ the corresponding highest weight, one gets $Rv_j^\pm=\pm v_j^\pm$, $Hv_j^\pm=jv_j^\pm$ and $F^+v_j^\pm=E^+v_j^\pm=0$.
The dimension of the representation  $[j]^\pm$ is $4j+1$ and the value of the Casimir $Q$ \eqref{eq:Q} is $\pm \frac{4j+1}{8}$.

In particular, in the fundamental bosonic representation $\left[\frac{1}{2}\right]^+$, the generators of $\mathfrak{osp}(1|2)$ are represented by the 
following $3\times3$ matrices
\begin{eqnarray}
 H=\frac{1}{2}\begin{pmatrix}
1 &0 &0\\
0 &-1&0\\
 0&0&0    
   \end{pmatrix}\quad,\quad  
   F^+=\frac{1}{2}\begin{pmatrix}0&0 &1\\
 0&0&0\\
 0&1&0    
   \end{pmatrix}\quad,\quad  
   F^-=\frac{1}{2}\begin{pmatrix}0&0 &0\\
 0&0&-1\\
 1&0&0    
   \end{pmatrix}\ ,
\end{eqnarray}
and $E^\pm=\pm 4(F^\pm)^2$, $R=\text{diag}(1,1,-1)$. 
For the sake of simplicity, we have used the same notations for the abstract algebraic elements of $\mathfrak{osp}(1|2)$ and their representatives. 

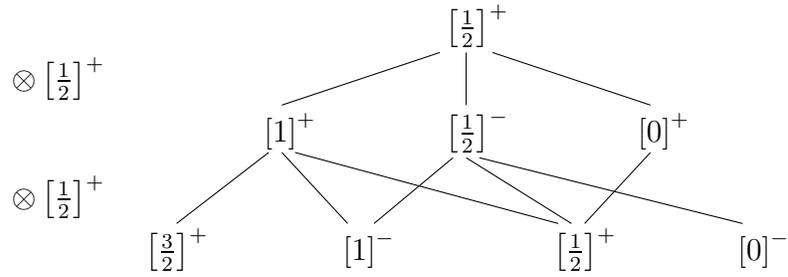
\begin{figure}[htb]
\begin{center}
\begin{tikzpicture}[scale=0.35]
\node at (0.9,-0.5) {$\left[\frac{1}{2}\right]^+$};
\draw (-0.5,-1.5) -- (-6.5,-3.5);\draw (0.5,-1.5)--(0.5,-3.5);\draw (1.5,-1.5) -- (7.5,-3.5);
\node at (-6.2,-4.5) {$\left[1\right]^+$};
\node at (0.9,-4.5) {$\left[\frac{1}{2}\right]^-$};
\node at (8,-4.5) {$\left[0\right]^+$};

\draw (-7,-5.3) -- (-10.5,-8); \draw (-6.5,-5.3) -- (-4,-8);  \draw (-6,-5.3) -- (4,-8);   

\draw (0,-5.5) -- (-3,-8);\draw (0.5,-5.5) -- (4.5,-8);\draw (1,-5.5) -- (11,-8);

\draw (7.5,-5.3) -- (5,-8);

\node at (-10.5,-9) {$\left[\frac{3}{2}\right]^+$};
\node at (-3.2,-9) {$\left[1\right]^-$};
\node at (5,-9) {$\left[\frac{1}{2}\right]^+$};
\node at (11.8,-9) {$\left[0\right]^-$};

\node at (-15,-2.5) {$\otimes\left[\frac{1}{2}\right]^+$};
\node at (-15,-7) {$\otimes\left[\frac{1}{2}\right]^+$};
\end{tikzpicture}
\caption{Bratteli diagram for the threefold tensor product of the fundamental representation. 
\label{fig:br}  }
\end{center}
\end{figure}

The direct sum decomposition of the tensor product of representations is also well-known \cite{PA1975,SNR1977,ERSS}. 
For the purpose of this paper, we need the following:
\begin{eqnarray}
 &&\left[0\right]^+ \otimes \left[\frac{1}{2}\right]^+ =  \left[\frac{1}{2}\right]^+\ ,\\
 &&\left[\frac{1}{2}\right]^+ \otimes \left[\frac{1}{2}\right]^\pm = \left[1\right]^\pm \oplus \left[\frac{1}{2}\right]^\mp\oplus  \left[0\right]^\pm\ ,\label{eq:dec2}\\
 &&\left[1\right]^+ \otimes \left[\frac{1}{2}\right]^+ = \left[\frac{3}{2}\right]^+ \oplus \left[1\right]^-\oplus  \left[\frac{1}{2}\right]^+\ .
\end{eqnarray}
With this information, we can draw the Bratteli diagram (see Figure \ref{fig:br}) that represents the direct sum decomposition of the threefold tensor product.

\vspace*{-0.5\baselineskip}
From this Bratteli diagram, we observe that 
\begin{equation}
\left[\frac{1}{2}\right]^+ \otimes \left[\frac{1}{2}\right]^+\otimes \left[\frac{1}{2}\right]^+ = 
  \left[\frac{3}{2}\right]^+ \oplus 2\left[1\right]^-\oplus  3\left[\frac{1}{2}\right]^+\oplus  \left[0\right]^-\ .\label{eq:dec4}\\
\end{equation}
We recall that the degeneracy of a representation present in the direct sum decomposition 
is the number of edges reaching the representation in the Bratteli diagram.   

\subsection{Centralizer of the threefold tensor product of the fundamental representation \label{sec:centr}}

Let us introduce $\mathcal{V}=\left[\frac{1}{2}\right]^+\otimes \left[\frac{1}{2}\right]^+ \otimes \left[\frac{1}{2}\right]^+$ and 
the centralizer associated to the action of $\mathfrak{osp}(1|2)$ on $\mathcal{V}$:
\begin{eqnarray}
 \mathcal{C}&=&\mathrm{End}_{\,\mathfrak{osp}(1|2)}\left(\mathcal{V}\right)\\
 &=&\left\{ M \in \mathrm{End}\left(\mathcal{V}\right)\ | \ 
 M.(g.v)=g.(M.v) \ , \ \forall v\in\mathcal{V}\ , \ \forall g\in \mathfrak{osp}(1|2) \right\}\ .
\end{eqnarray}
By adding the squares of the degeneracies ($1,2,3,1$) in the decomposition \eqref{eq:dec4}, one observes that the dimension of $\mathcal{C}$ is 15.
It is also known that the centralizer $\mathcal{C}$ is isomorphic to the Brauer algebra. Therefore, from Theorem \ref{thm:iso},
$\mathcal{C}$ is isomorphic to the quotiented Bannai--Ito algebra:
\begin{corollary}\label{co:cent}
The quotiented Bannai--Ito algebra $\overline{I}_3$ defined in Theorem \ref{thm:iso} 
is isomorphic to $\mathrm{End}_{\,\mathfrak{osp}(1|2)}\left(\left[\frac{1}{2}\right]^+\otimes \left[\frac{1}{2}\right]^+ \otimes \left[\frac{1}{2}\right]^+\right)$.
\end{corollary}
We can prove this corollary directly without reference to the isomorphism between the Brauer algebra and the centralizer $\mathcal{C}$.
In the following, we use the same notation, namely,  $Q_\#$, $X$, $Y$, $Z$ and $\omega$, for the algebraic elements and 
their images in $\mathrm{End}\left(\mathcal{V}\right)$.
Proposition \ref{pr:Phi} and relation \eqref{eq:centu2} imply that $X$, $Y$, $Z$ and $\omega$ are in $\mathcal{C}$. 
We must also show that the images in $\mathrm{End}\left(\mathcal{V}\right)$ of the l.h.s. of relations \eqref{eq:q1}--\eqref{eq:q5} vanish.
The Casimirs $Q_i$ (for $i=1,2,3$) equal $\frac{3}{8}$ times the identity matrix. 
From the decomposition of the tensor product of two fundamental representations 
into a direct sum of irreducible representations (see relation \eqref{eq:dec2}), 
we deduce that the eigenvalues of $Q_{12}$, $Q_{13}$ and $Q_{23}$ are $\frac{5}{8},-\frac{3}{8},\frac{1}{8}$. 
Therefore, from Proposition \ref{pr:Phi}, the eigenvalues of $X$, $Y$ and $Z$ are $-2,2,0$.
By the Cayley--Hamilton theorem, we conclude that relation \eqref{eq:q1} holds. 
We find similarly that the eigenvalues of $Q_4$ are given by $-\frac{5}{8},-\frac{1}{8},\frac{3}{8},\frac{7}{8}$ 
(see the third row of the Bratteli diagram displayed on Figure \ref{fig:br}) and that the eigenvalues of $\omega=\frac{7}{2}+12Q_4$ are $-4,2,8,14$. This proves
relation \eqref{eq:q2}. The eigenvalues of $X-\omega$ are given by the edges of the Bratteli diagram Fig.\ \ref{fig:br}: if $x$ is an eigenvalue of $X$ 
corresponding to the representation $[j]^{\epsilon_1}$ in the second row of the Bratteli diagram and $w$ is an eigenvalue of $\omega$ 
associated to the representation $[k]^{\epsilon_2}$ in the third row, then $x-w$ is an eigenvalue of $X-\omega$ iff $[j]^{\epsilon_1}$ and $[k]^{\epsilon_2}$ are 
connected in the Bratteli diagram. It is found this way that the eigenvalues of $X-\omega$ are $-16,-10,-8,-6,0,2,6$. This proves \eqref{eq:q3}. 
Relations \eqref{eq:q4} and \eqref{eq:q5} are derived similarly. This shows that the map from the quotiented Bannai--Ito algebra to $\mathcal{C}$ is 
a well-defined algebra homomorphism. 
The images of the following 15 elements $1$, $Q_{12}$, $Q_{23}$, $Q_{12}^2$, $Q_{23}^2$, $Q_{12}Q_{23}$, $Q_{23}Q_{12}$, $Q_{12}^2Q_{23}$, $Q_{12}Q_{23}Q_{12}$, $Q_{23}Q_{12}^2$, $Q_{23}^2Q_{12}$, 
$Q_{23}Q_{12}Q_{23}$, $Q_{12}^2Q_{23}^2$, $Q_{23}^2Q_{12}^2$ and $Q_{12}Q_{23}^2Q_{12}$ in $\mathrm{End}\left(\mathcal{V}\right)$ are linearly independent.
Surjectivity is therefore ensured since the dimension of the centralizer is 15.
In the proof of Theorem \ref{thm:iso}, we also show that $\mathrm{dim}(\overline{I}_3)=15$ which proves bijectivity.

\section{Conjectures and perspectives \label{sec:conj}}

Corollary \ref{co:cent} provides a link between a quotient of the Bannai--Ito algebra and the centralizer of the tensor product of three fundamental representations of $\mathfrak{osp}(1|2)$.
We believe that such a relation also exists for three copies of $\mathfrak{osp}(1|2)$ in the irreducible representation $[j]^+$. 
More precisely, we state the following conjecture:
\begin{conj}
Let $[j]^+$ be the irreducible representation of $\mathfrak{osp}(1|2)$ with $2j \in \mathbb{Z}_{>0}$.
The centralizer $\mathrm{End}_{\,\mathfrak{osp}(1|2)}\left([j]^+\otimes [j]^+ \otimes [j]^+ \right)$ is isomorphic to
the Bannai--Ito algebra $I_3$ defined by \eqref{eq:BI} quotiented by the following relations $\omega_X=\omega_Y=\omega_Z=\omega$ and
\begin{align}
& \prod_{k=-2j}^{2j}\big(X-2k\big)=0\,, \qquad \prod_{k=-2j}^{2j}\big(Y-2k\big)=0 \,, \qquad  \prod_{k=-2j}^{2j}\big(Z-2k\big)=0 \,, \label{eq:conj1}\\
& \prod_{k=-3j}^{3j} \big( \omega - (4j+1)(2j+1-2k) + 1 \big) = 0 \,, \label{eq:conj2} \\ 
& \prod_{k\in \mathcal{M}} \big(X-\omega -k\big)=0\,, \qquad 
\prod_{k\in \mathcal{M}} \big(Y-\omega -k\big)=0\,, \qquad 
\prod_{k\in \mathcal{M}} \big(Z-\omega -k\big)=0\,. \qquad   \,. \label{eq:conj3}
\end{align}
In the above formulas, the products are always understood to be with integer steps even if the boundaries have half-integer values.
The set $\mathcal{M}$ is obtained as explained at the end of the previous section from the edges of the 
Bratteli diagram associated to $[j]^+\otimes [j]^+ \otimes [j]^+$ (see below for an explicit description).
The isomorphism is defined by $(\pi_j\otimes \pi_j\otimes \pi_j)\Phi$ where $\Phi$ is given by \eqref{eq:Phi1} and \eqref{eq:Phi2} and $\pi_j$ is the representation homomorphism
from $\mathfrak{osp}(1|2)$ to $\mathrm{End}\left([j]^+\right)$.
\end{conj}
As an illustration, we give the Bratteli diagram for the threefold tensor product of the $[1]^+$ representation:
\begin{figure}[htb]
\begin{center}
\begin{tikzpicture}[scale=0.35]
\node at (7.5,-0.5) {$\left[1\right]^+$};
\draw (5.6,-1.5) -- (-4.9,-3.5);\draw (6.6,-1.5) -- (1.1,-3.5);\draw (7.1,-1.5)--(7.1,-3.5);\draw (7.6,-1.5) -- (13.1,-3.5);
\draw (8.6,-1.5) -- (19.1,-3.5);
\node at (-4.7,-4.5) {$\left[2\right]^+$};
\node at (1.4,-4.5) {$\left[\frac{3}{2}\right]^-$};
\node at (7.5,-4.5) {$\left[1\right]^+$};
\node at (13.6,-4.5) {$\left[\frac{1}{2}\right]^-$};
\node at (19.7,-4.5) {$\left[0\right]^+$};

\draw (-6.6,-5.5) -- (-17.1,-8);\draw (-5.6,-5.5) -- (-11.3,-8);\draw (-5.1,-5.5)--(-5.1,-8);\draw (-4.6,-5.5) -- (0.9,-8);\draw (-3.9,-5.5) -- (6.6,-8); 

\draw (-0.3,-5.5) -- (-10.7,-8);\draw (0.5,-5.5) -- (-5,-8);\draw (1,-5.5)--(1,-8);\draw (1.5,-5.5) -- (7,-8);\draw (2.2,-5.5) -- (12.8,-8); 

\draw (5.9,-5.5) -- (-4.6,-8);\draw (6.6,-5.5) -- (1.1,-8);\draw (7.1,-5.5)--(7.1,-8);\draw (7.6,-5.5) -- (13.1,-8);\draw (8.3,-5.5) -- (18.8,-8); 

\draw (12,-5.5) -- (1.5,-8);\draw (12.7,-5.5) -- (7.2,-8);\draw (13.2,-5.5)--(13.2,-8); 

\draw (19,-5.5) -- (7.7,-8); 

\node at (-16.9,-9) {$\left[3\right]^+$};
\node at (-10.8,-9) {$\left[\frac{5}{2}\right]^-$};
\node at (-4.7,-9) {$\left[2\right]^+$};
\node at (1.4,-9) {$\left[\frac{3}{2}\right]^-$};
\node at (7.5,-9) {$\left[1\right]^+$};
\node at (13.6,-9) {$\left[\frac{1}{2}\right]^-$};
\node at (19.7,-9) {$\left[0\right]^+$};

\node at (-24,-2.5) {$\otimes\left[1\right]^+$};
\node at (-24,-7) {$\otimes\left[1\right]^+$};
\end{tikzpicture}
\caption{Bratteli diagram for the threefold tensor product of the representation $[1]^+$.}
\end{center}
\end{figure}
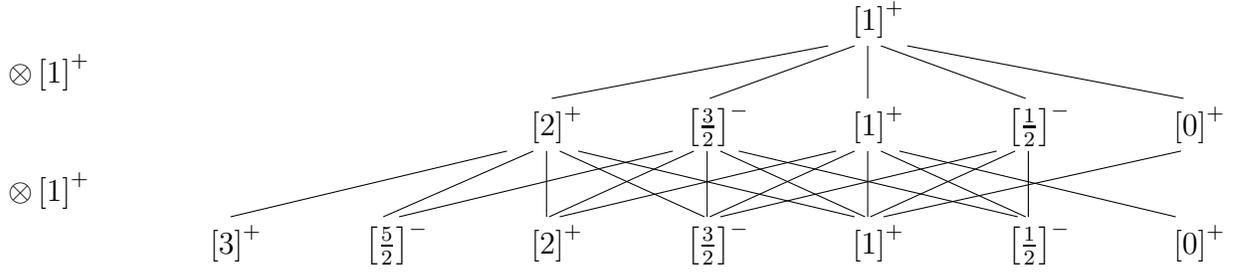

The eigenvalues of $X-\omega$ being given by the edges of the Bratteli diagram, one obtains for a generic $[j]^+$:   
\begin{align}
-6\ell-2k-5 \quad \text{with} &\quad |j-k| \le \ell \le j+k \;\; \text{and} \;\; 0 \le k \le 2j \label{eq:bradeb} \\
6\ell-2k-5 \quad \text{with} &\quad |j-k|+1 \le \ell \le j+k \;\; \text{and} \;\; 0 \le k \le 2j \\
\intertext{and, when $j$ is integer,}
-6\ell+2k-3 \quad \text{with} &\quad j-k \le \ell \le j+k \;\; \text{and} \;\; 0 \le k \le j-1 \nonumber \\
 \text{or} &\quad k-j+1 \le \ell \le j+k \;\; \text{and} \;\; j \le k \le 2j-1 \\
6\ell+2k-3 \quad \text{with} &\quad j-k \le \ell \le j+k+1 \;\; \text{and} \;\; 0 \le k \le j-1 \nonumber \\
 \text{or} &\quad k-j+1 \le \ell \le j+k+1 \;\; \text{and} \;\; j \le k \le 2j-1 
\intertext{or, when $j$ is half-integer,}
-6\ell+2k-3 \quad \text{with} &\quad j-k \le \ell \le j+k \;\; \text{and} \;\; 0 \le k \le j-\half \nonumber \\
 \text{or} &\quad k-j+1 \le \ell \le j+k \;\; \text{and} \;\; j+\half \le k \le 2j-1 \\
6\ell+2k-3 \quad \text{with} &\quad j-k \le \ell \le j+k+1 \;\; \text{and} \;\; 0 \le k \le j-\half \nonumber \\
 \text{or} &\quad k-j+1 \le \ell \le j+k+1 \;\; \text{and} \;\; j+\half \le k \le 2j-1.  \label{eq:brafin}
\end{align}
The total number of the $X-\omega$ eigenvalues, taking into account their multiplicities, is given by the Hex numbers $12j^2+6j+1$ (crystal ball sequence for hexagonal lattices). 
The set $\mathcal{M}$ is then obtained by considering the distinct eigenvalues given by equations \eqref{eq:bradeb}--\eqref{eq:brafin}. \\
The dimension of the centralizer $\mathrm{End}_{\,\mathfrak{osp}(1|2)}\left([j]^+\otimes [j]^+ \otimes [j]^+ \right)$ is equal to $d_j=(2j+1)^4-(2j)^4$, which is the sequence of rhombic dodecahedral numbers.

To support this conjecture, remark that for $j=\frac{1}{2}$ we recover the quotient of the Bannai--Ito algebra introduced in Theorem \ref{thm:iso}.
Another important point is that a similar conjecture has been made in \cite{CPV} for the connection between the Racah algebra 
and the centralizer of $\mathfrak{su}(2)$. In this case, the conjecture has been proven in numerous instances.
The main step to derive the conjectured isomorphism is to show that \eqref{eq:conj1}--\eqref{eq:conj3} generate the whole kernel.

If true, this conjecture gives a description of the centralizer for three copies 
of $\mathfrak{osp}(1|2)$ in the representation $[j]^+$.
We also believe that this conjecture can be generalized to the case of three arbitrary irreducible $\mathfrak{osp}(1|2)$ representations of finite dimension. 
It would be also interesting to consider tensor products of degree higher than three; this would connect to the
higher rank Bannai--Ito algebra that has been introduced in \cite{DGV} and comparisons could then be made with the limit $q \to 1$ of the algebra studied in \cite{LZ}.
Obviously, the generalization to the case of the quantum superalgebra should also be possible and a connection between the $q$-Bannai--Ito algebra \cite{GVZ} and 
the Birman--Murakami--Wenzl algebra \cite{BMW} is to be expected.\\

\textbf{Acknowledgements.} N.Cramp\'e is partially supported by Agence Nationale de la Recherche Projet AHA ANR-18-CE40-0001 and is gratefully holding a CRM--Simons professorship. 
L.Frappat warmly thanks the Centre de Recherches Math\'ematiques (CRM) for hospitality and support during his visit to Montreal in the course of this investigation. 
The research of L.Vinet is supported in part by a Natural Science and Engineering Council (NSERC) of Canada discovery grant.

\end{document}